\documentclass[preprint]{article}
\usepackage{amsthm}
 \usepackage{amssymb}
 \usepackage{amsmath}
 \usepackage{xcolor}
 \usepackage{soul} 
\newcommand{\nn}{\textbf n}
\newcommand{\Q}[0]{\mathbf{Q}}
\newcommand{\EE}{\mathbb E}
\newcommand{\e}[0]{\mathbf{e}}
\newcommand{\R}{\mathbb R}
\newcommand{\Z}[0]{\mathbb Z}

\newtheorem{thm}{Theorem}[section]
\newtheorem{prop}{Proposition}[section]
\newtheorem{lem}{Lemma}[section]
\newtheorem{cor}{Corollary}[section]
\newtheorem{defn}{Definition}[section]
\newtheorem{ex}{Example}[section]
\newtheorem{rem}{Remark}[section]

\title{Computable bounds on the spectral gap  for \\ unreliable Jackson networks} 

\author{Pawe{\l} Lorek \thanks{Work of both authors supported by NCN Research Grant DEC-2011/01/B/ST1/01305.
Postal address of both authors: Mathematical Institute, University of Wroc{\l}aw, pl. Grunwaldzki 2/4, 50-384
Wroc{\l}aw, Poland }
\thanks{E-mail: \texttt{Pawel.Lorek@math.uni.wroc.pl}}
\\{\it  University of Wroc{\l}aw} 
\and Ryszard Szekli \footnotemark[1] \thanks{E-mail: \texttt{Ryszard.Szekli@math.uni.wroc.pl}}
  \\{\it  University of Wroc{\l}aw}
}

\date{}
 
\begin{document}
\maketitle
\begin{abstract}
The goal of this paper is to identify  exponential convergence
rates and to find computable bounds for them for Markov  processes representing unreliable Jackson networks.
First we use the bounds of Lawler and Sokal in order to show that, for unreliable  Jackson networks,  the spectral gap is strictly
positive if and only if the spectral gaps for the corresponding coordinate birth and death processes are positive. Next, utilizing some  results on birth and death processes, we find  bounds on the spectral gap for network processes in terms of the hazard and equilibrium functions of the one dimensional
 marginal distributions of the stationary distribution of the network. These distributions must be  in this case strongly light-tailed, in the sense that  their discrete hazard functions have to be separated from zero. We relate these hazard functions with the corresponding networks' service rate functions using the equilibrium rates of the stationary one dimensional marginal distributions. We compare the obtained bounds on the spectral gap with some other known bounds. 
\end{abstract}

\section{Introduction}
We start with a description of the general setting used in this paper. Let ${\bf X}=( X_t, t\geq 0)$ be a Markov process on a countable state
 space $\EE$ with a bounded generator $\Q$ and the corresponding semi-group of operators $(P_t,\  t>0)$ on $L^2(\EE , \pi)$. We assume ergodicity of this process and the existence 
of the invariant probability measure $\pi$. 
The  usual scalar product on $L^2:=L^2(\EE , \pi)$ and the corresponding $L^2$ norm we denote by
$$\langle f,g\rangle _\pi = \sum_{\nn\in\EE} f(\nn)g(\nn)\pi(\nn),\qquad ||f||^2=\langle f,f\rangle _\pi ,$$
and by ${\pmb 1}$ the constant function equal to 1 on $\EE$. We shall use the symbol
$\pi(f)$ to denote $\langle f,\pmb{1}\rangle _\pi=E_\pi (f( X_t))$.
We denote the $L_2$ spectral gap corresponding to ${\bf X}$ by
\begin{equation}\label{deff:spectral_gap}
Gap(\Q) := \inf\left\{ -\langle f,\Q f\rangle _{\pi}: ||f||=1, \pi (f)=0\right\}.
\end{equation}
We say that ${\bf X}=( X_t, t\geq 0)$ has  an "exponential rate of convergence" if  $Gap(\Q)>0$. 

Then, for reversible processes,  
the following  conditions are equivalent (see, e.g., Theorem 1.9, \cite{Chen_ks_eigen})

 \begin{itemize}
 \item[(i)] 
 for all $f\in L^2(\EE , \pi ),$
 \begin{equation*}|| P_t f - \pi (f) || \le e^{-Gap(\Q) t} ||f-\pi (f)||,\  t>0,\end{equation*}
 \item[(ii)]
  for each $\e\in  \EE$ there exists $C(\e)>0$ such that $$||\delta_\e  P_t-\pi||_{tv} \leq C(\e) e^{-\alpha t },\  t>0, \ \text{for some}\  \alpha>0,$$
 \end{itemize}
where $||\cdot||_{tv}$ denotes the total variation norm. 

Denote by $\alpha_0$ the best rate in $||\delta_\e  P_t-\pi||_{tv}$ convergence. It is known that for ergodic birth and death processes $Gap(\Q)=\alpha_0$. See, e.g., \cite{vandoorn-81} or Theorem 5.3 in \cite{chen-l2-91}. We shall point out (section \ref{Speed of convergence to stationarity})  that we have this equality also for ergodic reversible (unreliable) Jackson networks. 

It is usually a very difficult (if not impossible) task to compute $Gap(\Q)$.
 Sometimes it is possible to prove that $Gap(\Q)>0)$ (the existence) without being able to give  computable bounds on the gap. We consider the problem of finding computable
 bounds for the $L_2$ spectral gap of unreliable Jackson network Markov processes which we will define later by the corresponding generators.

There exist very large literature on the speed of convergence to stationarity for general processes ${\bf X}$. Let us recall a few references. In order to prove the 
existence of the spectral gap for ${\bf X}$ it is possible to use the theory of Harris recurrent Markov processes, utilizing Lyapunov functions with appropriate drift conditions, see Meyn and Tweedie \cite{Meyn-Tweedie-93}.  However, computable bounds are not easily obtainable by the Harris recurrence techniques.  Some exceptions are known such as for example when $\EE=\R$ (totally ordered state space) and  in addition when the process is stochastically monotone, see \cite{LunMeyTwe}, \cite{Roberts-Tweedie-00}. Other approaches are possible via coupling methods or renewal theory methods. See, e.g., \cite{Chen_ks_eigen} , \cite{Baxendale-05}, \cite{Berenhaut-99},\cite{Berenhaut-02}. Sharper results leading to bounds on the spectral gap are possible via strong stationary times, strong stationary duality, Cheeger type inequalities,  Poincare 
inequalities or direct spectral representations for the semi-group $(P_t,\  t>0)$. See, e.g.,  \cite{LawSok88}, \cite{Ligget89},  \cite{Diaconis-Fill-90},
 \cite{Diaconis-Fill-90-aap}, \cite{Fill-91-peis}, \cite{Fill-92}, \cite{Fill-91}, \cite{Diaconis-Stroock-91}, \cite{KroSchTay}, \cite{LorSze12}, and in a book form see \cite{Chen_ks_eigen}.  
Symmetry assumptions turned out to be especially effective in analysis, and reversibility of ${\bf X}$ is a typical assumption for many results. However, even for
 birth and death processes analysis of  spectra and the transient behaviour of $(P_t, t>0)$ is far from being simple.  See, e.g.,  \cite{callaert-73},\cite{kijima-92}, \cite{kartashov-00}, 
\cite{vandoorn-02}, \cite{hart-03},\cite{sirl-07}, \cite{liu-09}, \cite{vandoorn-10}, \cite{chen-10}, \cite{chafai-12},  for some  results on bounds on the gap, and \cite{Diaconis-Miclo-09}, 
\cite{Fill-09}, \cite{LorSze12} for strong stationary times and duals  approach to finite state birth and death processes.
 
Jackson network processes can be seen as a generalization of birth and death processes, and one can expect that bounds for the spectral gap of a network should be related to some bounds on 
spectral gaps for some related birth and death processes. In fact,  Jackson network processes are much more complicated than  birth and death processes because they  are built upon an additional 
Markov chain which guides the routing inside the network. Reversibility for Jackson networks depends upon reversibility of the routing matrix. It is known that the simplest Jackson networks with 
constant service rates are stochastically monotone (under coordinate-wise ordering) but in general the stochastic monotonicity depends on the properties of the corresponding state dependent service
 rates.  See, e.g., \cite{daduna-szekli-95} for many monotonicity properties of Jackson networks. Unfortunately, for unreliable Jackson networks no reasonable stochastic monotonicity is present (see, e.g.,
\cite{DadKulSauSze06}),  therefore known methods to find computable bounds on the spectral gap, using the stochastic monotonicity property, are not applicable for networks (also because all known 
results on computable bounds with a use of stochastic monotonicity require totally ordered state spaces). A plausible expectation is that the speed of convergence to stationarity of a network 
should correspond to a bottleneck node of the network. Some partial results in this direction can be found for networks with state independent service rates in \cite{Anantharam-89} (for finite 
capacity networks), \cite{blanc-85}, and \cite{DieWarr} (for tandems). For networks with state independent rates also Lyapunov drift functions were studied in \cite{FayMal91}, \cite{HorSpi-92}.

A realted line of research is to study the essential spectrum of the generator $\Q$ of $(P_t,\  t>0)$. A broad view on this topic can be found in \cite{wu}. $\Q$ can act as operator on various function spaces (Banach lattices) such as for example $L^p, \  p\ge 1$, and the corresponding essential spectral gap is always larger than the gap defined by the underlying norm in a given function space. The essential spectral radius is directly related to LDP theory, to Lyapunov functions and asymptotic results for the tail distributions of (the first) returning times to compact sets. Finding the essential spectral radius for $L^2$ space gives at once an upper bound on the speed of convergence in $L^2$, which is interesting but more interesting for assessing the speed of convergence is to have lower bounds on the gap. In general, we do not know results characterizing when the $L^2$ spectral radius is equal to the corresponding essential spectral radius, however some examples showing this equality for some ergodic birth and death processes are known. See e.g., Example 8.4 in \cite{wu}.  For ergodic birth and death processes with constant (state independent) rates the $L^2$ spectral gap is known. See, e.g., examples after Corollary 1.3 in  \cite{chen-96}.  For ergodic birth and death processes with constant rates the essential $L^2$ spectral gap is also known. See, e.g.,  in the context of Jackson networks,  \cite{ignatiouk-06}, \cite{ignatiouk-12}. In the language of queueing processes, for an ergodic single $M/M/1$ station the $L^2$ spectral gap and the $L^2$ essential spectral gap are both equal to $(\sqrt\lambda-\sqrt \mu)^2$. It would be interesting to characterize the class of networks for which this equality holds true more generally. The fact that the problem of using spectral theory to characterise rates of convergence is a rather complex problem, even for countable Markov chains used in queueing theory, can be seen for example from \cite{malyshev-spieksma}, \cite{vere-jones} or \cite{Chen_ks_eigen}.
 
Positive lower bounds for the spectral gap of Jackson networks with state dependent service rates  were obtained  via some  related birth and death processes  in  \cite{McDon-94b}, by using 
conductance bounds from \cite{LawSok88}. 
A related comparison result for spectral gaps for classical Jackson networks is given in \cite{DadSze08}, Proposition 3.6, where a direct comparison involving the spectral gaps for some 
related birth and death processes is given, using an additional assumption on the routing. 
In this paper we give some bounds on the spectral gap for networks with state dependent service rates using Cheeger type constants using \cite{LawSok88}, similarly as 
in \cite{McDon-94b}, but related to some other birth and death processes than those defined in \cite{McDon-94b}. We consider in addition the possibility of  having unreliable 
nodes. Unreliable Jackson networks are networks, where in some subsets of the set of nodes the service stations
can be broken and then repaired during the time evolution of the system. The breakdown
and repair events can be of a rather general nature, but driven by a Markov process. In the time intervals  when nodes are
broken, there are several rules for re-routing.
 For full details of such networks see Sauer and  Daduna \cite{SauDad03}, and   Sauer \cite{SauerThesis}. We assume for unreliable networks reversibility, however this assumption  
can be skipped if the nodes are reliable. In a few  examples  we compare our bounds with bounds obtainable  from the results of  \cite{DadSze08} (lower bounds), and \cite{ignatiouk-12} (upper bounds).
Jackson networks possess two remarkable properties crucial for our analysis, namely the stationary distribution has a product form (also for unreliable networks) and exponential 
ergodicity for them is directly related to the strong light-tailness of the stationary distribution. It is worth mentioning that admitting service rates which are state dependent
 in the model implies that each discrete distribution with the support $\{0,1,2,\ldots\}$ can appear as the stationary distribution for a node in the network. We will characterize light-tailness of the stationary distribution by the corresponding discrete hazard rate functions. The stationary distribution can be also characterized by the corresponding so called equilibrium rates which turn out to be equal to individual, state dependent  traffic intensity functions for the nodes of a network.
Roughly speaking, the speed of convergence for a network will depend on a joint effect of how heavy the tails of the marginals of the stationary distribution are, together with how 
fast  each single node operates which in turn depends on the routing in the network.

The paper is organized as follows. In the next section we introduce unreliable networks by giving the respective generator. In section 3 we give a result relating the existence of the spectral gap of unreliable networks with the tail properties of its stationary distribution.  In section 4 we use equilibrium rates to reformulate our results from section 3. In section 5 we give the  proofs of the results from section 3. Finally, in section 6 we give some examples of bounds on the spectral gap for  networks.

\section{Description of the network process}\label{sec:def_process}
The classical {\bf Jackson
network} consists of $m$ numbered servers, denoted by
$ {M}:=\{1,\ldots,m\}$. Station $j\in {M}$ is a single server queue with
infinite waiting room under FCFS (First Come First Served) discipline. All the
customers in the network are indistinguishable. There is an external Poisson
arrival stream with intensity $\lambda$ and arriving customers are sent
to node $j$ with probability $r_{0j}$, $\sum_{j=1}^mr_{0j}=r\leq 1$. Customers arriving at node $j$ from the outside
or from other nodes request a service which is  at node $j$  provided with intensity $\mu_j(n)$
($\mu_j(0):=0$), where $n$ is the number of customers at node $j$ including the
one being served. All  service times and arrival processes are assumed to be
independent. \par A customer departing from node $i$ immediately proceeds to
node $j$ with probability $r_{ij}\geq 0$ or departs from the network with
probability $r_{i0}$. The routing is independent of the past of the system given the
momentary node where the customer is. We
assume that the stochastic  matrix $R:=(r_{ij}, \ i,j\in  {M}\cup\{0\})$ is irreducible.

Let $Z_j(t)$ be the number of customers present at node $j$, at time $t\geq 0$.
Then $$Z(t)=(Z_1(t),\ldots,Z_m(t))$$ is the joint queue length vector at time
instant $t\geq 0$ and ${\bf Z}:=(Z(t),t\geq 0)$ is the joint queue length
process with the state space $\EE=\mathbb{Z}_+^m$.

The unique stationary distribution for ${\bf Z}$ exists if and only if the unique solution of
the {\bf traffic equation}
\begin{equation}\label{eq:traffic}
\lambda_i=\lambda r_{0i}+\sum_{j=1}^m \lambda_j r_{ji}, \quad i=1,\ldots,m
\end{equation}
satisfies
$$C_i:= 1+\sum_{n=1}^\infty {{\lambda}_i^n\over \prod_{y=1}^n \mu_i(y)} <\infty, \quad 1\leq i \leq m.$$

The parameters of a Jackson network are: the arrival intensity $\lambda$,
the routing matrix $R$ (with the corresponding traffic arrival intensities vector ${\pmb
\lambda}=(\lambda_1,\ldots,\lambda_m)$), the vector of service rates ${\pmb
\mu}=(\mu_1(\cdot),\ldots,\mu_m(\cdot))$ and the number of servers $m$. Our standing assumption for all considered networks is that for all $j$, $\underline\mu_j:=\inf_{n\ge 1} \mu_j(n)>0$. We denote the overall minimal service intensity by $\underline\mu=\min_j \underline\mu_j$.

Assume now that the servers at the nodes in the Jackson network are
unreliable, i.e., the nodes may break down. The breakdown event may occur in
different ways. Nodes may break down as an isolated event or in groups
simultaneously, and the repair of the nodes may end for each node individually
or in groups as well. It is not required that those nodes which stopped service
simultaneously return to service at the same time instant. To describe the
system's evolution we have to enlarge the state space for the network process
as it will be described below. Denote by $ {M}_0:=\{0,1,\ldots,m\}$  the set of nodes enlarged by adding the \emph{outside} node.

Let $ {D}\subseteq  {M}$ be the set of servers out of order, i.e. in \emph{down status}.
\begin{itemize}
 \item  if $ {I}\subseteq {M}\setminus {D},   {I}\neq\emptyset$ is a subset of nodes in \emph{up status}, then all servers in $ {I}$ break down simultaneously with intensity
$\alpha^{ {D}}_{ {D}\cup {I}}( n_i: i\in  {M})$,

\item  
if $ {H}\subseteq {D},  {H}\neq\emptyset$, then all servers from
$ {H}$ return from repair simultaneously with intensity
$\beta_{ {D}\setminus {H}}^{ {D}}(n_i : i\in  {M})$.

\item The routing is changed according to so-called {\sc Repetitive Service -
Random Destination Blocking} (RS-RD BLOCKING) rule: For  $ {D}$  - set of
servers under repair routing probabilities are restricted to nodes from
$ {M}_0\setminus  {D}$ as follows:
$$
r^{ {D}}_{ij}=\left\{
\begin{array}{lll}
 r_{ij}, & i,j\in {M}_0\setminus  {D}, & i\neq j, \\
 r_{ii}+\sum_{k\in  {D}} r_{ik}, & i\in  {M}_0\setminus {D},& i=j.\\
\end{array}\right.
$$
The external arrival rates are
\begin{equation*}\label{external_arrivals}
\lambda r^{ {D}}_{0j}=\lambda r_{0j}
\mathrm{\  for\ nodes\  } j\in {M}\setminus {D},
\end{equation*}
and zero, otherwise. 
\end{itemize}
Let $R^{ {D}}=(r^{ {D}}_{ij})_{i,j\in {M}_0\setminus  {D}}$ be the modified routing.
Note that $R^{\emptyset}=R$.

 We assume for the intensities of breakdowns and repairs
 $\emptyset\neq  {I}\subseteq{ M\setminus {D}}$ and $\emptyset\neq {H}\subseteq  {D}$
 that
$$\begin{array}{lcr}
\alpha^{ {D}}_{ {D}\cup {I}}( n_i: i\in  {M}) &:=& {\psi( {D}\cup {I})\over \psi( {D})},\\
  & &\\
\beta^{ {D}}_{ {D}\setminus {H}}( n_i: i\in  {M}) &:=& {\phi( {D})\over \phi( {D}\setminus{H})},
\end{array}
$$
where $\psi$ and $\phi$ are arbitrary positive functions, defined for all subsets of the set of nodes, and $\psi(\emptyset)=\phi(\emptyset)=1$. That means that
 breakdown and
repair intensities depend on the sets of servers but are independent of the particular numbers of
customers present in these servers.

In order to  describe unreliable Jackson networks we need to attach to the
state space $\mathbb{Z}_+^m$ of the corresponding standard network process an
additional component which includes information on the availability of the system.
 We consider new state space
$$\nn =( {D}, n_1,n_2,\ldots,n_m)\in \mathcal{P}( {M})\times \mathbb{Z}_+^m =: {\bf\EE},$$
where $\mathcal{P}( {M})$  denotes the powerset of $ {M}$. The first (zero) coordinate in $\nn $ we call the availability coordinate.

The set $ {D}$ is the set of servers in \emph{down status}. At node $i\in {D}$
there are $n_i$ customers waiting for server to be repaired. Denote  possible transitions by

\begin{equation}\label{deff:transfText}
\begin{array}{lll}
T_{ij}\nn & := & ( {D},n_1,\ldots,n_i-1,\ldots,n_j+1,\ldots, n_m), \\
T_{0 j}\nn & := & ( {D},n_1,\ldots,n_j+1,\ldots, n_m), \\
T_{i0 }\nn & := & ( {D},n_1,\ldots,n_i-1,\ldots, n_m), \\
T_{ {H}}\nn & := & ( {D}\setminus {H},n_1,\ldots, n_m), \\
T^{ {I}}\nn & := & ( {D}\cup {I},n_1,\ldots, n_m). \\
\end{array}
\end{equation}
\begin{defn}\label{unreliable-net}
The Markov process ${\bf X}=({\bf X}(t),t\geq 0)$ defined by the  infinitesimal
generator
\begin{equation}\label{eq:jack_zawodne_gen}
{\small
\begin{array}{llll}
\displaystyle {\bf\Q}f({\bf\nn})=&\displaystyle \sum_{j=1}^m [f(T_{0
j}{\bf\nn})-f({\bf\nn})]  {\lambda} r^{  D}_{0j} & + &
\displaystyle \sum_{i=1}^m\sum_{j=1}^m [f(T_{ij}{\bf\nn})-f({\bf\nn})]
\mu_i(n_i)r^{  D}_{ij} + \\[15pt]
& \displaystyle \sum_{\emptyset\neq {I}\subseteq {M\setminus D}} [f(T^{ {I}}{\nn})-f({\bf\nn})]{ \psi( {D}\cup {I})\over \psi( {D})}& +&\displaystyle
\sum_{ \emptyset\neq {H}\subseteq {D}} [f(T_{ {H}}{\bf\nn})-f({\bf\nn})]
{ \phi( {D})\over \phi( {D}\setminus{H})}+\\[15pt]
&\displaystyle \sum_{j=1}^m[f(T_{j0 }{\bf\nn})-f({\bf\nn})]\mu_j(n_j)
r^{  D}_{j0} &
\\

\end{array}
}
\end{equation}
is called {\bf unreliable Jackson network}.
\end{defn}
We denote the corresponding transition intensities (written in a matrix form) by
$[q({\bf\nn},{\bf\nn}')]_{{\bf\nn},{\bf\nn}'\in\EE \ }.$

Similarly to the classical case the invariant distribution for this Markov process can be written in a product form.
\begin{thm}[Sauer and Daduna \cite{SauDad03}]\label{twr:unreliable_cornelia}
 Let ${\bf\bf X}$ be unreliable Jackson network  following
 the RS-RD-BLOCKING. If the routing matrix $R$
is reversible, i.e.:
$$\lambda_j r_{ji}=\lambda_ir_{ij}, \qquad i,j\in {M},$$
then the stationary distribution of process ${\bf\bf X}$ is given by

\begin{equation}\label{eq:stat_distr_unrel}
\pi({\bf\nn})=\pi( {D},n_1,\ldots,n_m)= {1\over C}{\psi( {D})\over
\phi( {D})} \prod_{i=1}^m \pi_i(n_i),
\end{equation}
 where
\begin{equation}\label{pi}\pi_i(n_i)={1\over C_i} {\lambda_i^{n_i} \over
\prod_{k=1}^{n_i}\mu_i(k)}, \qquad C_i=1+\sum_{n=1}^\infty {\lambda_i^n\over
\prod_{y=1}^{n}\mu_i(y)}\end{equation} and $C$ is the normalization constant used for the availability coordinate.
Constants $C_i, i=1,\ldots,m$ are all finite if and only if the network is ergodic.
 \end{thm}
Note that in this generality, the reduced state vector to the number of customers alone, without the availability coordinate, does not form a Markov process. The model of unreliable network is an analogue of the classical Jackson network model but it can not be reduced to the classical one by adjusting  parameters of the availability coordinate since all configurations of down nodes are possible with positive probability under our assumptions.

\subsection{Equilibrium rate and hazard rate for stationary distribution}
For a non-negative  random variable $X\in\Z_+$, with probability function
$p(k)=P(X=k)$, such that for any  $k\in\Z_+$, $P(X=k)>0$,
the \emph{total hazard  function} $H_p$ is defined for all $x\ge0$ by
$$H_p(x)=-\log \bar{F}(x).$$
Further, the \emph{discrete hazard function} we define for natural arguments by
\begin{equation}\label{eq:hazard_rate}
h_p(k)={p(k)\over \bar F(k-1)},\ k\ge 0,
\end{equation}
where $\bar{F}(k)=P(X>k)$. 
Note that for such a variable, for natural arguments $k\ge 0$
\begin{equation}\label{ogon_hazard}
H_p(k)=-\log \prod_{j=0}^k \left(1-h_p(j)\right).
\end{equation}
and for arbitrary $x\geq 0$ we have
\begin{equation}\label{ogon_hazardcont}
H_p(x)=-\log \prod_{j=0}^{\lfloor x\rfloor} \left(1-h_p(j)\right)=\sum_{j=0}^{\lfloor x\rfloor}\log\left(\frac{1}{1-h_p(j)}\right),
\end{equation}
where   $\lfloor x\rfloor$ denotes the integer part of $x$.

\begin{defn}
We say that a discrete distribution $(p(k), k=0,1,\ldots)$ (or a discrete random variable $X$) is \emph{\bf strongly light-tailed} if there exists $\epsilon>0$ such that $\inf _{k\ge 0} h_p(k)>\epsilon$.
\end{defn}

The following lemma and example explain how {\em the strong light-tailness} and {\em the usual light-tailness} are related.
Recall the usual light-tailness. An arbitrary distribution  function $F$ with its support contained in $[0,\infty)$ is light-tailed if $\int _0^\infty e^{sx}dF(x)<\infty$ for some $s>0$. 

\begin{lem}
Consider  a random variable $X\in\Z_+$, with probability function
$p(k)=P(X=k)$, such that for any  $k\in\Z_+$, $P(X=k)>0$, and $p$ is strongly light-tailed. Then it is  light-tailed in the usual sense.
\end{lem}
\begin{proof} It is known (see e.g. Rolski {\em et al.} \cite{RolSchSchTeu}, Th. 2.3.1) that $$\liminf _{x\to\infty} -\frac{1}{x}\log(\bar F(x))>0$$ implies that $F$ is light-tailed.
Note that
$$\frac{H_p(x)}{x}\geq \frac{H_p(\lfloor x\rfloor)}{\lfloor x\rfloor+1},$$
for all $x\ge 0$, therefore
\begin{equation}\label{l-t}
\inf_n \frac{H_p(n)}{n+1}>0 \Rightarrow \liminf_{x\to\infty}\frac{H_p(x)}{x}>0.
\end{equation}
From the exponential light-tailness we have   for all $j$, $\log(\frac{1}{1-h_p(j)})>\log(\frac{1}{1-\epsilon}),$ and hence from  (\ref{ogon_hazardcont})
$$
\frac{H_p(n)}{n+1}>\log\left(\frac{1}{1-\epsilon}\right)>0,
$$
which from (\ref{l-t}) implies that $F$ is light-tailed.

\end{proof}

We give now a simple example in order to see that for discrete distributions strong light-tailness is a strictly  stronger notion than the usual light-tailness. 
(This example shows at the same time that there exists a birth and death process having its rate of convergence to stationarity not exponentially fast, but having its stationary distribution  light-tailed).
\begin{ex}\label{examplel-t}
\rm
Let us take as $p$ the distribution which corresponds to the hazard function $h_p$ given by $h_p(1)=1/2,$

\begin{displaymath}
h_p(k)=\left\{ \begin{array}{ccllc}
1/k & \textrm{if} & k=2n+1,&\  n\ge 1,\\
1/2  & \textrm{if} & k=2n,&\  n\ge 0.\\
\end{array}\right.
\end{displaymath}
This distribution is not  strongly light-tailed since $\inf_k h_p(k)=0$. 
However, for each natural $n$, $\lim_{n\to\infty}\frac{H_p(2n+1)}{2n+2}=\lim_{n\to\infty}\frac{H_p(2n)}{2n+1}=\log(2)/2>0$, and from (\ref{l-t}) we obtain that $p$ is light-tailed. 
\end{ex} 

For a non-negative  random variable $X\in\Z_+$, with probability function
$p(k)=P(X=k)$, such that for any  $k\in\{0,1,2,\ldots\}$, $P(X=k)>0$, we define the
\emph{equilibrium rate} function  for natural arguments by
\begin{displaymath}
e_p(k)=\left\{ \begin{array}{ccccc}
{p(k+1)\over p(k)} & \textrm{if} & k\geq 0, \\
0  & \textrm{if} & k<0. \\
\end{array}\right.
\end{displaymath}

Since under our assumptions the equilibrium rate function $(e_p(k), k\ge 0)$ uniquely determines the probability function $(p(k), k\ge 0)$, it is therefore 
possible to express strong light-tailness in terms of equilibrium rates. The following formulas connect hazard and equilibrium rate functions

\begin{equation}\label{eq:hazard_equil}
e_p(k)=\frac{ h_p(k+1)(1-h_p(k))}{h_p(k)},\quad k\ge 0
\end{equation}
and
\begin{equation}\label{eq:hazard_equil2}
h_p(k)={1\over 1+\sum_{j=k}^\infty  e_p(k)\cdots e_p(j)},\quad k\ge
0.
\end{equation}

It is worth mentioning that each discrete distribution with the support $\Z_+$ can appear as the stationary distribution for a birth and death process 
with constant birth rates and variable death rates. Strong  light-tailness of $\pi_i$ can be expressed in terms of the corresponding equilibrium rates,
 which in turn  are equal to the corresponding birth/death ratios. A precise formulation for a single birth and death process we give in the following lemma.
\begin{lem}\label{twr:gap_const_lambda}
Consider  $\{p(k)\}_{k\ge 0}$ an arbitrary probability function on $\Z_+$, such that
$p(k)>0, k\ge 0$, with the corresponding equilibrium rate  $e_p(k), \ k\ge 0$. Then for each  birth and death process ${\bf Z}$ with fixed $\lambda(k)\equiv \lambda
>0,\ k\ge 0,$ and  death rates defined by $$\frac{\lambda}{\mu(k+1)}= e_p(k), \ k\ge 0,$$
the stationary distribution of ${\bf Z}$ is equal to $p(k),\ k\ge 0$,
\end{lem}
\begin{proof} For the stationary distribution $\check\pi$ of the birth death process ${\bf Z}$  we have
$$\check\pi(i)/\check\pi(0)={\lambda^i\over \mu(1)\cdots \mu(n)}={\lambda^i\over \lambda^i {p(0)\over p(1)} {p(1)\over p(2)}\cdots {p(i-1)\over p(i)}} =
p(i)/p(0), \ i\ge 1.$$
Thus we have $p=\check\pi$.
\end{proof}

Neither $h_p(k)$ nor $e_p(k)$ have to be convergent as $k\to\infty$. However,
from (\ref{eq:hazard_equil}), (\ref{eq:hazard_equil2}) we obtain a connection between these limits if they exist and are finite.
\begin{lem}\label{arbitrary}
  Consider  $\{p(k)\}_{k\ge 0}$, an arbitrary probability function on $\Z_+$, such that
$p(k)>0, k\ge 0$, with the corresponding equilibrium rate  $e_p(k), \ k\ge 0$. Then

$h_p=\lim_{k\to\infty} h_p(k) $ exists and $h_p\in (0,1)$ if and only if $e_p=\lim_{k\to\infty} e_p(k)$ exists and $e_p\in (0,1).$
In this case
$$h_p=1- e_p.$$
\end{lem}
\begin{ex}
Recall that the negative binomial distribution is defined by
$$
p(k)={r+k-1\choose k}(1-p)^kp^r, \ \  r>0, \ \ k=0,1,\ldots , \ p\in (0,1)
$$
The corresponding equilibrium rate is given by 
$$
e_p(k)=(1-p)(k+r)/(k+1), \ \ k=0,1,\ldots .
$$
The corresponding limit at infinity fulfills $e_p=(1-p)$, and for the corresponding limit at infinity of the hazard rate we get $h_p=p>0$, 
which means that this distribution is strongly light-tailed.
\end{ex}
\begin{ex}
For the Poisson distribution
$$
p(k)=e^{-\lambda}\lambda^k/k!,\ \ \lambda>0,\ \ k=0,1,\ldots,
$$
and $$e_p(k)=\frac{\lambda}{k+1}.$$
For the corresponding limits at infinity we have here $e_p=0$, and $h_p=1$, the Poisson distribution is strongly light-tailed.
\end{ex}
It is worth mentioning that the negative binomial and Poisson distributions fit into the so called Panjer recurrence scheme, more precisely, 
we say that $p(k)$ fulfills Panjer's recurrence if for some $a,b\in \R$
$$
p(k+1)=\left(a+\frac{b}{k+1}\right)p(k),\ \  k=0,1,\ldots,
$$
which is equivalent to saying that the corresponding equilibrium rate has a  hyperbolic form 
$$
e_p(k)=a+\frac{b}{k+1}.
$$  
For the negative binomial distribution $a:=1-p$, and $b:=(r-1)(1-p)$.
In both cases the equilibrium rate function is monotone. Distributions with non-increasing equilibrium rates are equivalently called $PF_2$ densities,
 for more details in connection with queueing networks see \cite{daduna-szekli-96}. 
\begin{ex}
A discrete analog of the Pareto distribution can be defined by
$$
p(k)=C\frac{1}{(k+1)^\alpha}, \ \ \alpha>1, \ \ k=0,1,\ldots,
$$
where $C$ is the normalization constant. Then
$$
e_p(k)=\left ( \frac{k+1}{k+2}\right ) ^\alpha.
$$
For the corresponding limits at infinity we have here $e_p=1$, and $h_p=0$, this distribution is heavy-tailed.

\end{ex}
In the context of unreliable queueing networks it is natural to define the ratio $\frac{\lambda_i}{\mu_i(k+1)}$, being a function of $k$ variable,  as 
{\em the traffic intensity function} for the $i$-th station.
From lemma \ref{twr:gap_const_lambda} it follows that for ergodic networks the traffic 
intensity function at the $i$-th station is equal to the equilibrium rate of the marginal $\pi_i$ distribution of the network's stationary 
distribution $\pi$.  If we assume that the service intensity at node $i$ is non-decreasing as a function of the number of customers at
 this node, then $\pi_i$ has a $PF_2$ density, and it is strongly light-tailed. Another possibility is that the traffic intensity function is
 increasing to 1 at a selected node $i$, and the network is ergodic but having at the node $i$ a heavy-tailed distribution $\pi_i$. 
It will be showed in the next section that in such a case the network process will not converge to stationarity geometrically fast. 
Also,  if at a fixed station $i$ the traffic intensity function is not monotone and corresponds to a light-tailed distribution which is not 
strongly light-tailed as in example  \ref{examplel-t}, then such a network also will not converge to stationarity 
geometrically fast.

\section{Existence of spectral gap and light tailed distributions}
\begin{thm}\label{twr:jack_unrel_gap_main}
\hfill

\noindent (i)  Let ${\bf X}$ be ergodic unreliable Jackson network process following
 the RS-RD- BLOCKING, with the infinitesimal generator ${\Q}$. Suppose that ${\Q}$ is bounded and the minimal service intensity $\underline\mu>0$. 

If the routing matrix $R$
is reversible   then  $Gap({\Q})>0$ if and only if  all distributions $\pi_i, \ i=1,\ldots,m$ are strongly light-tailed.

\noindent (ii)  Let ${\bf Z}$ be ergodic classical Jackson network process
  with the corresponding infinitesimal generator ${\Q(\bf Z)}$. Suppose that ${\Q(\bf Z)}$ is bounded and the minimal service intensity $\underline\mu>0$. 

Then  $Gap({\Q(\bf Z)})>0$ if and only if  all distributions $\pi_i, \ i=1,\ldots,m$ are strongly light-tailed.

\end{thm}

The proof of this theorem will be given in  section \ref{proof}.

We formulated the results on the positivity of the spectral gap and on the convergence to stationarity in terms of the discrete hazard functions of
 the stationary distribution. For queueing networks it would be however more reasonable to formulate the assumptions in terms of the parameters of the network. 

The existence of the spectral gap of an unreliable network can be formulated in terms of the corresponding arrival and service rates (as a consequence of  Theorem \ref{twr:jack_unrel_gap_main} and  Lemma \ref{twr:gap_const_lambda}) as follows 
\begin{cor}\label{cor-rates}
Let ${\bf X}$ be an ergodic unreliable Jackson network process following
 the RS-RD-BLOCKING, with the infinitesimal generator ${\Q}$. Suppose that ${\Q}$ is bounded and the minimal service intensity $\underline\mu>0$. If the routing matrix $R$
is reversible then  $Gap({\Q})>0$ if and only if for each $i=1,\ldots,m,$ 
$$\inf_k{1\over 1+\sum_{j=k+1}^\infty {\lambda_i^{j-k} \over \mu_i(k+1)\cdots \mu_i(j)}}>0.$$
In particular for ergodic networks, if for all $i=1,\ldots,m,$  the limits for the traffic intensity  functions $\lim_{k\to\infty}\lambda_i/\mu_i(k)<1$ exist then $Gap({\Q})>0.$
\end{cor}
 For the classical Jackson networks the assumption on reversibility can be skipped.
\subsection{Speed of convergence to stationarity}\label{Speed of convergence to stationarity}

Denote by $\alpha_0$ the best rate in $||\delta_\e  P_t-\pi||_{tv}$ convergence. It is known that for ergodic birth and death processes $Gap(\Q)=\alpha_0$, see e.g. \cite{vandoorn-81} or Theorem 5.3 in \cite{chen-l2-91}. From Theorem 8.8. (2) \cite{Chen_ks_eigen},  for ergodic reversible processes it is known that $\alpha_0\ge Gap(Q)$.  From Theorem 8.13, (4) \cite{Chen_ks_eigen}, we have

\begin{thm}\label{twr:jack_unrel_gap_main_rate}
 Let ${\bf X}$ be an ergodic, unreliable Jackson network following
 the RS-RD-BLOCKING, with generator ${\Q}$,
 given by (\ref{eq:jack_zawodne_gen}), and the corresponding transition semigroup
 $( P_t)$. Suppose the routing matrix $R$
is reversible.

If  $\pi_i$ is strongly light-tailed, for each $ i=1,\ldots, m$, then the following conditions are equivalent
\begin{itemize}

 \item[(i)] for all $f\in L^2(\EE,\pi)$
$$|| P_tf - \pi(f) || \leq e^{-Gap({\Q}) t} ||f-\pi(f)||,\ t>0,$$

\item[(ii)] for each $\e\in  \EE$ there exists $C(\e)>0$ such that $$||\delta_\e  P_t-\pi||_{tv} \leq C(\e)
e^{-Gap({\Q})t },\ t>0,$$
\end{itemize}
where $||\cdot||_{tv}$ denotes the total variation norm.
\end{thm}
\begin{proof}
First note that the network process is reversible under the assumption that $R$ is reversible. It is enough to check the assumptions of Theorem 8.13, (4) \cite{Chen_ks_eigen}. Let  $p_t(\e,\e')=\frac{dP_t(\e,\cdot)}{d\pi}(\e'), \ t>0, \e,\e'\in  \EE$. Then $p_{2s}(\e,\e) =P({\bf X}(2s)=\e|{\bf X}(0)=\e)/\pi(\e)$. Hence $p_{2s}(\cdot,\cdot)\in L^{(1/2)}_{loc}(\pi)$ (with the usual notation for $L^p(\pi)$ spaces as in \cite{Chen_ks_eigen}) if $\sum_{\e\in A\subset\EE}(\pi(\e))^{(1/2)}<\infty$ for bounded $A$, which trivially holds. The set of bounded functions with compact support is (also trivially) dense in $L^2(\pi)$ since $\EE$ is a discrete space.

\end{proof}
\begin{rem} For the classical Jackson networks, the reversibility assumption on the routing matrix $R$ can be relaxed in order to obtain the implication $(i)\Rightarrow (ii)$.
\end{rem}

\section{Bounds on the spectral gap}
In this section we recall some bounds on the spectral gaps of birth and death processes.  
For a more complete description  see \cite{Chen_ks_eigen}, (chapter 5), \cite{chen-10}, \cite{vandoorn-02}, \cite{vandoorn-10}, and references therein.

Let us recall Theorem 3.7 of 
Liggett \cite{Ligget89}. For convenience we give formulation of it  simplified to the case of state independent birth rates. 
\begin{thm}[Liggett  \cite{Ligget89}]\label{twr:ligett_thr_3.7}
 Assume that ${\bf Z}$ is an ergodic birth and death process on $\Z_+$,  with state independent birth rates $\lambda >0$, and possibly state dependent death rates $\mu(n)>0$, and
 for all $i\geq 0$,  and for some $c, d>0$, we have
$$\sum_{j>i} \pi(j)\leq c\pi(i)\lambda  \quad \mathrm{and} \quad \sum_{j>i}\pi(j)\leq d\pi(i).$$
Then for the corresponding generator $\Q({\bf Z})$,
\begin{equation}\label{ligget-bound}
Gap(\Q({\bf Z}))\geq { (\sqrt{d+1}-\sqrt{d})^2\over c } \geq {1\over 2c(1+2d)}.
\end{equation}
\end{thm}
In the case of constant birth rates, from the Corollary 3.8 of Liggett  \cite{Ligget89}, we have that a necessary and sufficient condition for
$Gap(\Q({\bf Z}))$ to be positive is that the stationary distribution is such that $$\inf_{i\ge 0} \frac{\pi(i)}{\sum _{j\ge i}\pi (j)}>0,$$ 
which is by definition the strong light-tailness. 
\noindent 
Therefore from  Corollary 3.8 of Liggett  \cite{Ligget89} we have
\begin{lem}\label{b-d-bound}
Assume that ${\bf Z}$ is an ergodic birth and death process on $\Z_+$,  with state independent birth rates $\lambda >0$, and possibly state dependent death rates $\mu(n)>0$. Then $Gap(\Q({\bf Z}))>0$ if and only if the stationary distribution $\pi$ is strongly light tailed. Moreover, if for  some $\epsilon>0$
 we have
$$\inf _{n\ge 0}h_{\pi}(n)\ge \epsilon, $$
then 
\begin{equation}\label{b-d-bound-f}
Gap(\Q({\bf Z}))\geq \frac{\lambda (1-\sqrt{1-\epsilon})^2}{1-\epsilon}\geq \frac{\lambda \epsilon^2}{2(1-\epsilon)(2-\epsilon)}.
\end{equation}
\end{lem}

\begin{proof} From $\sum_{j>i} \pi(j)\leq c\pi(i)\lambda $ we have $\sum_{j\ge i} \pi(j)\leq c\pi(i)\lambda +\pi(i)$, so for the lower bound on the hazard function we have $\epsilon=1/(1+c\lambda)$, therefore $c=(1-\epsilon)/(\lambda \epsilon)$. Similarly we get $d=(1-\epsilon)/\epsilon$, and using (\ref{ligget-bound}) we obtain (\ref{b-d-bound-f}).
\end{proof}

A lower bound on the spectral gap can be given directly in terms of the birth and death rates. See, e.g., \cite{vandoorn-02}.
\begin{lem}\label{vandoorn-bound}
Assume that ${\bf Z}$ is an ergodic birth and death process on $\Z_+$,  with state independent birth rates $\lambda >0$, and possibly state dependent death rates $\mu(n)>0$. Then
$$
Gap(\Q({\bf Z}))\ge \inf_{n\ge 0}\left [\lambda +\mu(n+1) -\sqrt{\lambda \mu(n)} - \sqrt{\lambda \mu(n+1)}\right ].
$$
\end{lem}
\begin{rem}
For more details on estimation of spectral gaps for birth and death processes see Corollary 1.2, Corollary 1.3 in \cite{chen-96}, and also \cite{chen-10}, \cite{Chen_ks_eigen}, \cite{vandoorn-85}, \cite{vandoorn-02}, \cite{vandoorn-10}. It is natural to ask how  do different bounds compare. It turns out that  optimality of a given bound strongly depends on the parameters of a given birth-death process, as described in an example after Theorem (5.2), \cite{ chen-l2-91}. In a sense, different bounds are {\sl incomparable} - as stated there. For particular cases it is reasonable to try out all possibilities.
\end{rem}
Combining the above bounds for birth and death processes and the bounds obtained in the proof of Theorem \ref{twr:jack_unrel_gap_main} (see (\ref{final-1})) we have from (\ref{b-d-bound-f}) 
\begin{prop}\label{final}
\hfill

\noindent (i) Let ${\bf X}$ be an ergodic, unreliable Jackson network following
 the RS-RD-BLOCKING, with generator ${\Q}$,
 given by (\ref{eq:jack_zawodne_gen}). Suppose the routing matrix $R$
is reversible.

If  $\pi_i$ is strongly light-tailed, for each $ i=1,\ldots, m$, and  
$$\inf _{n\ge 0}h_{\pi_i}(n)\ge \epsilon_i>0, $$
then

\begin{equation}\label{final-1-1}
Gap({\Q}) \geq \frac{1}{8|{\Q}|}\left ({{{q}}^{min}\over {\check{q}}^{max}}\frac{Gap(\check\Q_0)\wedge\displaystyle\min_{1\le i\le m}\frac{\lambda_i (1-\sqrt{1-\epsilon_i})^2}{1-\epsilon_i} }{1+\bar d \hspace{1pt}  \bar b (2m +1) }\right )^2\ \ \ \ \  
\end{equation}
and
$$
\begin{array}{l}
Gap({\Q}) \geq \\
\frac{1}{8|{\Q}|}\displaystyle\left ({{{q}}^{min}\over {\check{q}}^{max}}\frac{Gap(\check\Q_0)\wedge\displaystyle\min_{1\le i\le m} \inf_{n\ge 0}\left [\lambda_i +\mu_i(n+1) -\sqrt{\lambda_i \mu_i(n)} - \sqrt{\lambda_i \mu_i(n+1)}\right ] }{1+\bar d \hspace{1pt}  \bar b (2m +1) }\right )^2,\\ 
\end{array}
$$
where $\bar d,\ \bar b,\ |{\Q}|, q^{min}, \check{q} ^{max}$ are defined by (\ref{factor-d}), (\ref{factor-b}), (\ref{factor-norm}), (\ref{factor-q}), (\ref{factor-q-check}), respectively.
\medskip

\noindent (ii)  Let ${\bf Z}$ be ergodic classical Jackson network process
  with the corresponding infinitesimal generator ${\Q(\bf Z)}$. Suppose that ${\Q(\bf Z)}$ is bounded and the minimal service intensity $\underline\mu>0$. 
If  $\pi_i$ is strongly light-tailed, for each $ i=1,\ldots, m$, and  
$$\inf _{n\ge 0}h_{\pi_i}(n)\ge \epsilon_i>0, $$
then 
\begin{equation}\label{final-1-2}
Gap({\Q(\bf Z)})\geq \frac{1}{8|{\Q(\bf Z)}|}\left ({{{q}}^{min}\over {\check{q}}^{max}}\frac{\displaystyle\min_{1\le i\le m}\frac{\lambda_i (1-\sqrt{1-\epsilon_i})^2}{1-\epsilon_i} }{1+ \bar b 2m  }\right )^2
\end{equation}
and
$$
\begin{array}{l}
Gap({\Q}) \geq  \\
\frac{1}{8|{\Q(\bf Z)}|}\displaystyle\left ({{{q}}^{min}\over {\check{q}}^{max}}\frac{\displaystyle\min_{1\le i\le m}\inf_{n\ge 0}\left [\lambda_i +\mu_i(n+1) -\sqrt{\lambda_i \mu_i(n)} - \sqrt{\lambda_i \mu_i(n+1)}\right ]}{1+ \bar b 2m  }\right )^2. \hfill\\ 
\end{array}
$$
%
\end{prop}
In all above given bounds the factor $1+\bar d \hspace{1pt} \bar b (2m +1)$ can be reduced to 1 if in the network $r_{i0}>0$ and $r_{0i}>0$ for all $i=1,\ldots, m.$
The bounds obtained in the above proposition are valid for a quite general class of networks but it is reasonable to search for alternative bounds and alternative methods under some additional structural assumptions. 
 We recall two cases for classical Jackson networks, the first one with state dependent service rates but fulfilling a partial balance requirement for the routing matrix (see \cite{DadSze08}, Proposition 4.4), the second one for classical Jackson networks with state independent service rates (see \cite{ignatiouk-12}).
\begin{prop}\label{gap2}
Let ${\bf Z}$ be ergodic classical Jackson network process
  with the corresponding infinitesimal generator ${\Q(\bf Z)}$. Suppose that ${\Q(\bf Z)}$ is bounded and the minimal service intensity $\underline\mu>0$. 
 Assume that the routing matrix $R$ has strict positive departure probabilities $r_{i0}>0$ and that $\lambda r_{0i}>0$ for  $i=1,\dots,m$. 

 Assume further a partial balance condition 
 \begin{equation}\label{detbal1}
 \lambda_j \sum_{i=1}^m  r_{ji}= \sum_{i=1}^m\lambda_i r_{ij},\quad \forall j=1,\dots,m.
 \end{equation}

Then 
$$ Gap(\Q({\bf Z}))  \geq \min_{1\le i \le m}Gap(\tilde Q_i) ,$$
where, for $i=1,\ldots, m$,  $\tilde Q_i$ denotes the generator of the birth and death process with the birth rate $\lambda r_{0i}$ and the state dependent death rate $\mu_i(n_i)r_{i0}$.
\end{prop}
\begin{cor}
Under the assumptions of Proposition \ref{gap2}, if in addition   $\pi_i$ is strongly light-tailed, for each $ i=1,\ldots, m$, and  
$$\inf _{n\ge 0}h_{\pi_i}(n)\ge \epsilon_i>0, $$
 then
\begin{equation}\label{daduna-1}
Gap({\Q({\bf Z})}) \geq \displaystyle\min_{1\le i\le m}\frac{\lambda r_{0i} (1-\sqrt{1-\epsilon_i})^2}{1-\epsilon_i} 
\end{equation}
and
\begin{equation*}
\displaystyle Gap({\Q({\bf Z})}) \geq \min_{1\le i\le m}\inf_{n\ge 0}\left [\lambda r_{0i}+\mu_i(n+1)r_{i0} -\sqrt{\lambda r_{0i} \mu_i(n)r_{i0} } - \sqrt{\lambda r_{0i}\mu_i(n+1)r_{i0}}\right ].\hfill
\end{equation*}
\end{cor}
Now we  recall from \cite{ignatiouk-12} some special cases of classical Jackson networks in order to present some (upper) bounds on the corresponding $L^2$ spectral gap. The results in \cite{ignatiouk-12} are related to the essential spectral gap. We shall compare  in section \ref{numerical} our lower bounds with the presented below  upper bounds and will obtain in some cases a nice approximation for $L^2$ spectral gap.
Because the essential $L^2$ spectral gap is larger then  $L^2$ spectral gap we have from Corollary 3.4, and Proposition 3.6 in \cite{ignatiouk-12}.
\begin{prop}\label{gap3}
Let ${\bf Z}$ be ergodic classical Jackson network process
  with the corresponding infinitesimal generator ${\Q(\bf Z)}$. Assume that the service intensities are state independent.

\noindent (i) If the routing is completely symmetrical, i.e. $r_{ij}=p<1/(m-1)$ for all $i\neq j$, $i,j=1,\ldots, m$,
and for some $i_0 \in \{1,\ldots, m\}$ we have
\begin{equation}\label{one}
\min_{1\le i \le m}(\sqrt{\mu_i}-\sqrt{\lambda_i})=\sqrt{\mu_{i_0}}-\sqrt{\lambda_{i_0}}
\end{equation}
and
\begin{equation}\label{two}
\min_{1\le i \le m}\left (  \frac{\mu_i}{\sqrt{\mu_{i_0}}}- \frac{\lambda_i}{\sqrt{\lambda_{i_0}}}     \right)=\sqrt{\mu_{i_0}}-\sqrt{\lambda_{i_0}},
\end{equation}
then
$$
Gap({\Q(\bf Z)})\le \left ( 1-\frac{(m-1)p^2}{1-(m-2)p}\right )\min_{1\le i \le m}(\sqrt{\mu_i}-\sqrt{\lambda_i})^2.
$$

\noindent (ii) If $m=3$, and 
\begin{equation}\label{eq:matrixRm3}
 R=\left (\begin{array}{cccc}
0 & r_{01} & r_{02} & r_{03} \\
1-(p+q)&0&p&q\\
1-(p+q)&q&0&p\\
1-(p+q)&p&q&0
\end{array}\right ),
\end{equation}
where $p,q\in (0,1)$, $p+q<1$, then

$$
Gap({\Q({\bf Z})})\le \frac{1-p^3-q^3-3pq}{1-pq}\min_{1\le i \le m} (\sqrt{\mu_i}-\sqrt{\lambda_i})^2 
$$
provided $\lambda_i/\mu_i=\lambda_j/\mu_j,\  i,j\in M$  or there exists $i_0$  such that $\mu_i\ge \mu_{i_0}$	and $\lambda_i \le  \lambda_{i_0}$, for all $i$.

\end{prop}

\section{Proof of  Theorem \ref{twr:jack_unrel_gap_main}}\label{proof}
We give the proof of Theorem \ref{twr:jack_unrel_gap_main} using the following theorem.
\begin{thm}[Liggett  \cite{Ligget89}, Th. 2.6]\label{liggett:twr2.6}
Suppose that a pure jump Markov process ${\bf X}$, with generator $\check\Q$ and
stationary distribution $\pi$ evolves on the product state space
$\EE=\EE_0\times \EE_1\times\cdots \EE_m$, $m\ge 1$, having coordinates which are
independent Markov processes such that $i-$th coordinate has generator $\check\Q_i$,
denumerable state space $\EE_i$ and  invariant probability measure $\pi_i$.
Then $\pi$ is the product measure of $\pi_i$'s and
$$Gap(\check\Q)=\min_{0\le i\le m} Gap(\check\Q_i).$$
\end{thm} 

Proof of Theorem \ref{twr:jack_unrel_gap_main} (i). We assume that the availability coordinate process is not degenerate with $\phi$ and $\psi$  positive.
Let $\check{\Q}$ be the generator associated with $(m+1)$-dimensional process
$(\mathbf{Y}_t, \check{\mathbf{Z}}_t)_{t\ge 0}$, where $\check{\mathbf{Z}}_t$ is the vector of
$m$ independent birth and death processes with generators $ \check{\Q}_i$,
$i=1,\ldots,m$, given by
\begin{equation}\label{eq:Qhati}
\check{\Q}_i f(n)=[f(n+1)-f(n)]\lambda_i+[f(n)-f(n-1)]\mu_i(n), \quad n\in 
\mathbb{N},
\end{equation}
and  $\mathbf{Y}_t$ is the process on state space $\mathcal{P}( {M})$ with
infinitesimal generator denoted by $\check\Q_0$ and the stationary distribution:
$$\pi_{0}( {I})={1\over C} { \psi( {I})\over \phi( {I})},
\qquad C:= \left(\sum_{ {I}\subseteq {M}}{\psi( {I})\over
\phi( {I})}\right).$$

We write $[\check{q}(\nn,\nn')]_{\nn,\nn'\in {\EE}}$ for the corresponding transition intensities.

The stationary distribution of the process with generator $\check{\Q}_i$  is $\pi_i$, which is given in the product formula (\ref{pi}) for networks.

Consider the following Cheeger's constants for $A\subset \EE$

$${\kappa}(A):={\sum_{\nn\in A} \pi(\nn){{q}}(\nn,A^c)\over
\pi(A)\pi(A^c)}, \qquad {\kappa}:=\inf_{A:\pi(A)\in(0,1)} {\kappa}(A),$$

$$\check{\kappa}(A):={\sum_{\nn\in A} \pi(\nn){\check{q}}(\nn,A^c)\over
\pi(A)\pi(A^c)}, \qquad \check{\kappa}:=\inf_{A: \pi(A)\in(0,1)} \check{\kappa}(A), $$
where $\pi$ is given by  (\ref{eq:stat_distr_unrel}).

We will show that there exist $0<v_1, v_2<\infty$ such that uniformly for all $A\subset \EE$
\begin{equation}\label{eq:v1v2_exist}
  v_2 \sum_{\nn \in A}\pi(\nn)\check{q}(\nn,A^c)\geq \sum_{\nn\in A}\pi(\nn){q}(\nn,A^c)\geq v_1 \sum_{\nn\in A}\pi(\nn)\check{q}(\nn,A^c).
\end{equation}

 Then with $0<v_1, v_2<\infty$ as in (\ref{eq:v1v2_exist}), we use Theorem 2.1 in Lawler and Sokal  \cite{LawSok88}, and since the process with the generator  $\check{\mathcal{\Q}}$ is
 reversible, we have that $Gap(\check{\mathcal{\Q}})\le \check{\kappa}$. Further,
 uniformly in  $A$, $\check{\kappa}(A)\le (v_1)^{-1}
{\kappa}(A)$, hence $\check{\kappa}\le (v_1)^{-1}
{\kappa}$. Under our assumptions we will have $Gap(\check{\mathcal{\Q}})>0$ which in turn, using Theorem 2.3 in Lawler and Sokal  \cite{LawSok88} (which assures
 that $\kappa^2/(8|{\Q}|)\le Gap({\Q}) $)  will imply  that $Gap({\Q})>0$. Here 
\begin{equation}\label{factor-norm}
|{\Q}|=\pi-{\rm ess\  sup}_{\nn}\  q(\nn, \{\nn\}^c).
\end{equation}
Similarly, it is possible to argue that $Gap({\Q})>0$ implies that $Gap(\check{\mathcal{\Q}})>0$.

In order to complete the proof we turn now to show the validity of (\ref{eq:v1v2_exist}) which is equivalent to
\begin{equation}\label{eq:Jackson_unrel_infsup}
 \inf_{A\subset \EE\atop \pi(A)\in(0,1)} \left\{
{\sum_{\nn\in A}\pi(\nn){q}(\nn,A^c) \over
 \sum_{\nn\in A}\pi(\nn)\check{q}(\nn,A^c)
}\right\}\ge v_1>0 
\end{equation}
and
\begin{equation}\label{equation-sup}
 \sup_{A\subset \EE\atop \pi(A)\in(0,1)} \left\{
{\sum_{\nn\in A}\pi(\nn){q}(\nn,A^c) \over
 \sum_{\nn\in A}\pi(\nn)\check{q}(\nn,A^c)
}\right\}\le v_2<\infty.
\end{equation}

For a fixed $A$, such that $\pi(A)\in (0,1)$, we define
$$\partial {A}=\{\nn\in A: {{q}}(\nn,A^c)>0\}, \qquad \partial \check{A}=\{\nn\in A: {\check{q}}(\nn,A^c)>0\}.$$

Let
\begin{equation}\label{factor-q}
{{q}}^{min}=\inf_{A:\pi(A)\in (0,1)}\inf_{\nn\in \partial{A} } \left\{{{q}}(\nn,A^c)\right\}, \quad
  {{q}}^{max}=\sup_{A:\pi(A)\in (0,1)}\sup_{\nn \in \partial{A}}\left\{{{q}}(\nn,A^c)\right\}.
\end{equation}
From our assumptions  the generators are bounded and $\underline\mu>0,$ therefore
$ {{q}}^{min}>0,$ and
 ${{q}}^{max}<\infty$.

For
\begin{equation}\label{factor-q-check}
{\check{q}}^{min}=\inf_{A:\pi(A)\in (0,1)}\inf_{\nn\in \partial\check{A}}  \left\{{\check{q}}(\nn ,A^c)\right\}, \quad
  {\check{q}}^{max}=\sup_{A:\pi(A)\in(0,1)}\sup_{\nn\in \partial\check{A} }\left\{{\check{q}}(\nn ,A^c)\right\},
\end{equation}
we also have ${\check{q}}^{min}>0$ and ${\check{q}}^{max}<\infty.$

For each $A$ such that $\pi(A)\in(0,1)$, we have
$$
{\sum_{\nn\in A}\pi(\nn){q}(\nn,A^c) \over
 \sum_{\nn\in A}\pi(\nn)\check{q}(\nn,A^c)
}=
{\sum_{\nn\in
\partial{A}}\pi(\nn){{q}}(\nn,A^c) \over 
\sum_{\nn\in \partial\check{A}} \pi(\nn){\check{q}}(\nn, A^c)
 },$$
so we obtain
$$
{{{q}}^{max}\over {\check{q}}^{min}}\cdot { \sum_{\nn\in
\partial{{A}}}\pi(\nn) \over  \sum_{\nn\in
\partial\check{A}} \pi(\nn)} \geq { \sum_{\nn\in
\partial{{A}}}\pi(\nn){{q}}(\nn, A^c) \over  \sum_{\nn\in
\partial\check{A}} \pi(\nn){\check{q}}(\nn,A^c)
 }
\geq {{{q}}^{min}\over {\check{q}}^{max}} \cdot {
\sum_{\nn\in
\partial{{A}}}\pi(\nn) \over  \sum_{\nn\in
\partial\check{A}} \pi(\nn)}.
$$
We shall continue our argument in the case of the lower bound (\ref{eq:Jackson_unrel_infsup}). The existence of this  lower bound ensures that if  $Gap(\check{\mathcal{\Q}})>0$,
 then   $Gap({\Q})>0$. Note that from Theorem \ref{liggett:twr2.6}, and Lemma \ref{b-d-bound}, the inequality $Gap(\check{\mathcal{\Q}})>0$ is equivalent to the condition that 
$\pi_i$ is strongly light-tailed, for each $ i=1,\ldots, m$. The proof for the upper bound is similar and we skip it.
In order to show (\ref{eq:Jackson_unrel_infsup}) it is enough to check that 
\begin{equation}\label{eq:Jackson_pi_partial}
 0< \inf_{A:\pi (A)\in (0,1)}\zeta(A)\ \  \quad
\textrm{where } \ \  \zeta(A):= { \sum_{\nn\in \partial{A}}\pi(\nn) \over  \sum_{\nn\in \partial{\check{A}}} \pi(\nn)}.
\end{equation}

If the network is such that for all $i=1,\ldots, m$, $r_{0i}>0$ and $r_{i0}
>0$ then $\partial\check{A}\subseteq \partial{A}$. In that case $\inf_{A:\pi (A)\in (0,1)}\zeta(A)\ge 1 $, and we can take $v_1={{{q}}^{min}\over {\check{q}}^{max}} $. Otherwise, 
we have to analyse $\partial\check{A}$, and $\partial{A}$ in more detail.

Let us examine the difference between $\pi(\nn)$ and
$\pi(\nn')$  when $\nn'$ and  $\nn$ differ exactly on one \textcolor{black}{nonavailability} coordinate by at most 1 and when $\nn$ and $\nn'$ have two different  sets of broken nodes $D$, and $D'$.

Recall from (\ref{eq:stat_distr_unrel}) that for $\nn=( {D},n_1,\ldots,n_m)\in
\mathcal{P}( {M})\times \mathbb{Z}_+^m$ we have:
$$
\pi(\nn)=\pi( {D},n_1,\ldots,n_m)= {1\over C}{\psi( {D})\over
\phi( {D})} \prod_{i=1}^m \pi_i(n_i),
 \qquad \mathrm{ \ where \ }
\pi_i(n_i):={1\over C_i} {\lambda_i^{n_i} \over \prod_{y=1}^{n_i}\mu_i(y)}.
$$
For  $n_i\geq 1,$
$$\pi_i(n_i+1)={1\over C_i} {\lambda_i^{n_i+1} \over \prod_{y=1}^{n_i+1}\mu_i(y)}
=\pi_i(n_i) {\lambda_i\over \mu_i(n_i+1)}$$
and
$$\pi_i(n_i-1)={1\over C_i} {\lambda_i^{n_i-1} \over {\prod_{y=1}^{n_i-1}\mu_i(y)}}
=\pi_i(n_i) {\mu_i(n_i)\over \lambda_i},
$$
thus, using $\underline\mu_i:=\inf_n \mu_i(n)>0$  and $\bar\mu_i:=\sup_n \mu_i(n)<\infty$, we have bounds  

$$
\begin{array}{rccll}
{\lambda_i\over \bar\mu_i}\pi_i(n_i) & \leq & \pi_i(n_i+1) & \leq
\pi_i(n_i){\lambda_i\over \underline\mu_i}, \\
\\
 {\underline\mu_i\over\lambda_i } \pi_i(n_i) & \leq &
\pi_i(n_i-1) & \leq \pi_i(n_i){\bar\mu_i\over \lambda_i}. \\
\end{array}
$$
Define
\begin{equation}\label{factor-b}
\bar b=\max_{1\le i\leq m}\left( \frac{\bar\mu_i}{\lambda_i
} \right),\  \underline b=\min_{1\leq i\leq m}\left( \frac{\lambda_i} {\bar\mu_i} \right),\   
\end{equation}

\begin{equation}\label{factor-d}
\bar d=\max_{D_1\neq D_2}\frac{\psi(D_2)\phi(D_1)}{\phi(D_2)\psi(D_1)}\ \  \text{and}\ \  \underline d=\min_{D_1\neq D_2}\frac{\psi(D_2)\phi(D_1)}{\phi(D_2)\psi(D_1)}
\end{equation}
Then,  if $\nn$ and $\nn'$  differ by at most $ 1$ on \textcolor{black}{exactly one}  
coordinate $i\in\{1,\ldots, m\}$, and have sets $D,\ D'$ on the availability coordinate then
\begin{equation}\label{eq:theta}
 \underline b\pi_i(n_i) \leq \pi_i(n_i')  \leq \bar b\pi_i(n_i)
\end{equation}
and  
\begin{equation}\label{eq:Theta}
 \underline d \hspace{1pt} \underline b\pi(\nn) \leq \pi (\nn')  \leq \bar d \hspace{1pt} \bar b \pi(\nn).
\end{equation}

We  rewrite $\zeta(A)$ as
$$
\zeta(A)= { \sum_{\nn\in
\partial{A}\cap\partial{\check{A}}}\pi(\nn) +\sum_{\nn\in
\partial{{A}}\setminus\partial\check{A}}\pi(\nn)  \over   \sum_{\nn\in
\partial\check{A}\cap\partial{{A}}}\pi(\nn) +\sum_{\nn\in
\partial\check{A}\setminus\partial{{A}}}\pi(\nn)   }.$$ 
 
Let us  consider  $\nn\in\partial{\check{A}}\setminus\partial{{A}}$.
Then there exists some $\nn ' \in A^c$ such that original process with the
intensity ${q}$ cannot move there in one step, but the process with 
$\check{q}$ can. The state $\nn '$ must be of the form $\nn '=T_{0i_0}\nn$ or $\nn '=T_{j_00}\nn$ (arrival or departure) \textcolor{black}{ since changing
availability coordinate is always possible in both processes, i.e., either both processes would leave $A$ or none.} 
We shall analyse 
the case of arrival since in the case of departure we can argue analogously. The key observation in this argument is the following:  if $\nn '=T_{0i_0}\nn$,  but the 
arrival intensity to node $i_0$ is equal to zero for the network process or this arrival movement is blocked by $D$ then the node $i_0$ must be reachable by an {\em unblocking} 
movement $D\to \emptyset$ and then $T_{0i_0}$ transition, or by an unblocking movement $D\to \emptyset$ and then an arrival to some station different than $i_0$, and a  migration
 movement or a series of consecutive migration movements. There are possibly multiple paths, but we can search for the minimal ones (which can  be multiple with the same length). 
Intuitively speaking we search for the 
shortest connection to a source node (i.e., a node which admits arrivals from the outside) from $i_0$ node (in the case of departure movement $\nn '=T_{j_00}\nn$ we search for 
the shortest connection to a sink node).   Consider all shortest paths of movements that connect $\nn$ with $\nn'$ in the network.
Denote such a path by $ \nn=\nn _0, \textcolor{black}{\nn_1=T_D\nn_0,}\ldots, \nn_k\textcolor{black}{=T^D\nn_{k-1}}=\nn ' \  \  (k\leq m\textcolor{black}{+1})$. Note that each such a path is not greater than $m+1$ since we can take as the first transition the 
one which puts $D$ to $\emptyset$ on the availability coordinate,  and  the worst case for the other transitions is when the station $i_0$ is the last station in a $m-$ series network.
Moreover,   each state on the path differs
 from $\nn$ by at most 1 on only one  non-availability coordinate (because on non-availability coordinates an arrival changes one coordinate by plus 1, and consecutive transitions
 change coordinates in such a way that after a transition the resulting state has exactly one coordinate changed by plus 1). 
Further, there  exists  a state $\nn _j$ on this path such that the network process  leaves $A$, and 
either  $\nn _j\in\partial\check{A}\cap\partial{{A}}$ or $\nn _j\in\partial{A}\setminus\partial\check{{A}}$.  Since $\nn_j$ differs from  $\nn$  by at most 1 on exactly one coordinate,
from (\ref{eq:Theta}) we have $\pi(\nn)\leq \bar d \bar b \pi(\nn_j)$. If we take two points on the border $\partial\check{A}\setminus\partial{{A}}$ for which the coordinate-wise 
distance is big enough, then the corresponding  border points  on $\partial{A}$ defined above  must be different, because $\nn_j$ always differs from  $\nn$  by at most 1 on a single 
coordinate. More precisely, let $\nn\in\partial{\check{A}}\setminus\partial{{A}}$ and ${\bf m}\in\partial{\check{A}}\setminus\partial{{A}}$ are such that they are different by more
 than two on each coordinate then the corresponding points $\nn_j$ and ${\bf m}_{j'}$, elements of $\partial{A}$, are distinct. In order to give a very rough bound 
on $\sum_{\partial{\check{A}}\setminus\partial{{A}}}\pi(\nn)$ we observe that for a fixed $\nn_j$ point there are not more than $2m+1$ points that are different 
by at most one on a single  coordinate from $\nn_j$, and $\nn_j$ can potentially be on a transition  (unblocking and migration) path described above for these 
points. Therefore we have
$$ \sum_{\partial{\check{A}}\setminus\partial{{A}}}\pi(\nn)\leq 
\bar d \hspace{1pt} \bar b (2m+1) \left(  \sum_{\nn\in \partial\check{A}\cap\partial{{A}}}\pi(\nn)
+\sum_{\nn\in
\partial{{A}}\setminus\partial\check{A}}\pi(\nn)\right) 
$$
and
$$
\begin{array}{lll}
\zeta(A) & \ge &{ 
\sum_{\nn\in \partial\check{A}\cap\partial{{A}}}\pi(\nn) +\sum_{\nn\in
\partial{{A}}\setminus\partial\check{A}}\pi(\nn)  \over 
\sum_{\nn\in \partial\check{A}\cap\partial{{A}}}\pi(\nn) +\bar d \hspace{1pt} \bar b (2m+1)
\left( \sum_{\nn\in \partial\check{A}\cap\partial{{A}}}\pi(\nn)
+\sum_{\nn\in
\partial{{A}}\setminus\partial\check{A}}\pi(\nn) \right)        }\\[24pt]
& \ge & { 
\sum_{\nn\in \partial\check{A}\cap\partial{{A}}}\pi(\nn) +\sum_{\nn\in
\partial{{A}}\setminus\partial\check{A}}\pi(\nn) \over (1+\bar d\hspace{1pt}  \bar b (2m+1) ) \left( \sum_{\nn\in
\partial\check{A}\cap\partial{{A}}}\pi(\nn) +\sum_{\nn\in
\partial{{A}}\setminus\partial\check{A}}\pi(\nn) \right)        } =\frac{1}{1+\bar d \hspace{1pt} \bar b (2m+1) }.\end{array}
$$
Summing up we obtain
$$
{ \sum_{\nn\in
\partial{{A}}}\pi(\nn){{q}}(\nn, A^c) \over  \sum_{\nn\in
\partial\check{A}} \pi(\nn){\check{q}}(\nn,A^c)
 }
\geq {{{q}}^{min}\over {\check{q}}^{max}} \cdot {
\sum_{\nn\in
\partial{A}}\pi(\nn) \over  \sum_{\nn\in
\partial\check{A}} \pi(\nn)}\ge {{{q}}^{min}\over {\check{q}}^{max}}\cdot \frac{1}{1+\bar d \hspace{1pt} \bar b (2m+1)} 
$$
and
 $$\check{\kappa}{{{q}}^{min}\over {\check{q}}^{max}}\cdot\frac{1}{1+\bar d \hspace{1pt} \bar b (2m+1) }\le 
{\kappa},$$
which implies (using Theorem 2.3 in Lawler and Sokal  \cite{LawSok88})
$$
 Gap({\Q})\geq \left (\check{\kappa}{{{q}}^{min}\over {\check{q}}^{max}}\cdot\frac{1}{1+\bar d \hspace{1pt} \bar b (2m+1)  }\right )^2/(8|{\Q}|) ,
$$
$$
 Gap({\Q})\geq \left ({{{q}}^{min}\over {\check{q}}^{max}}\cdot\frac{Gap(\check\Q)}{1+\bar d \hspace{1pt} \bar b (2m+1) }\right )^2/(8|{\Q}|)  
$$
and finally
\begin{equation}\label{final-1}
Gap({\Q}) \geq \left ({{{q}}^{min}\over {\check{q}}^{max}}\cdot\frac{\min_{0\le i\le m} Gap(\check\Q_i)}{1+\bar d \hspace{1pt} \bar b (2m+1)  }\right )^2/(8|{\Q}|).
\end{equation}
Proof of (ii). Note that we cannot specify parameters of an ergodic unreliable Jackson network process ${\mathbf X}$ to obtain the classical ergodic 
Jackson network process ${\mathbf Z}$ as a special case. However, it is possible to repeat all steps in the proof of (i)  for ${\mathbf Z}$ (skipping 
the availability coordinate, and reducing $2m+1$ to $2m$)  to get
\begin{equation}\label{final-2}
Gap(\Q({\bf Z}))\geq \left ({{{q}}^{min}\over {\check{q}}^{max}}\cdot\frac{\min_{1\le i\le m} Gap(\check\Q_i)}{1+\bar b 2m  }\right )^2/(8|\Q({\bf Z})|).
\end{equation}

\section{Numerical examples}\label{numerical}
We shall use two examples from \cite{ignatiouk-12} in order to to estimate $L^2$ spectral gap.
\begin{ex}
Let $\mathbf{Z}$ be the classical Jackson network with $m=3$ stations with  the  arrival intensity $\lambda$ and the routing matrix $R$ given in (\ref{eq:matrixRm3}), and  with $r_{01}=r_{02}=r_{03}=1/3$, where $p,q\in(0,1), p+q<1$.
Then $\lambda_1=\lambda_2=\lambda_3=\lambda/(3(1-(p+q))$ is the solution to the traffic equation.
Moreover, assume that service intensities are constant and are given by $\mu_i=c\lambda_i, i=1,2,3,$ where $c>1$. The network is ergodic
with stationary distribution being the product of $\pi_i , i=1,2,3$, where $\pi_i(k)=(1-\frac{1}{c})(\frac{1}{c})^{k}, i=1,2,3, \ k=0,1,\ldots$. The conditions
of Proposition \ref{gap3} $(ii)$ are fulfilled and we have:
$$Gap({\bf Q}({\bf Z})) \le Gap_{ess}:={1-p^3-q^3-3pq\over 1-pq} \lambda_1\left(\sqrt{c}-1\right)^2=$$$$ {p^2+p-pq+q^2+q+1\over  1-pq} {\lambda\over 3} (\sqrt{c}-1)^2.$$
We will compare the above upper bound  with the bounds given in Proposition \ref{gap2} and Proposition \ref{final}.
\medskip\par 
Let us start with the bound given in Proposition \ref{gap2}.
The partial balance condition (\ref{detbal1})  holds, and all birth and death processes $\tilde{Q}_i, i=1,2,3$ are equal in distribution. Denote the 
arrival intensity of  $\tilde{Q}_i$ process  by $\tilde{\lambda}_i$,  and its service rate by $\tilde{\mu_i}$.
We have  $\tilde{\lambda}_i=\lambda r_{0i} = \lambda/3$ and $\tilde{\mu}_i=\mu_i r_{i0}=c\lambda /3$.
As already indicated in the introduction the formula for $L^2$ spectral gap, for ergodic birth and death processes with constant rates, is known. The $L^2$ spectral gap (and the corresponding essential spectral gap) of $\tilde{Q}_i$ is given by
$$ Gap(\tilde{Q}_i)=\left(\sqrt{\tilde{\mu}_i}-\sqrt{\tilde{\lambda_i}}\right)^2={\lambda\over 3}(\sqrt{c}-1)^2,$$
therefore the resulting bound is
$$Gap({\bf Q}({\bf Z})) \geq {\lambda\over 3}(\sqrt{c}-1)^2.$$
It is worth mentioning that this bound does not depend on $p,q$. Moreover,
$$\inf_{p,q\in(0,1)\atop p+q<1} Gap_{ess}:={\lambda\over 3}(\sqrt{c}-1)^2.$$
 On the other hand,
$$\sup_{p,q\in(0,1)\atop p+q<1} Gap_{ess}= \lambda (\sqrt{c}-1)^2,$$
which means that the bound given in Proposition \ref{gap2} is at most 3 times smaller than the considered upper bound on the spectral gap. Moreover, the spectral gap $Gap({\bf Q}({\bf Z})) $ is arbitrarily close to ${\lambda\over 3}(\sqrt{c}-1)^2$ for small values of $p$ and $q$.

\medskip\par Now, let us turn to Proposition \ref{final}.
Each distribution $\pi_i$ is  geometric with the corresponding hazard functions  $h_{\pi_i}(n)=1-\frac{1}{c}$.
We have $r_{i0}>0$ and $r_{0i}>0$ for $i=1,2,3.$, thus we can reduce $1+\bar{d}\hspace{1pt}\bar{b}(2m+1)$ to 1 in this proposition.
We need yet to calculate:
$$
\begin{array}{lll}
|{\bf Q}| & =  & \lambda r_{01}+\lambda r_{02}+\lambda r_{03}+\mu_1+\mu_2+\mu_3=\lambda + 3c{\lambda\over 3(1-(p+q)}=\lambda\left(1+{c\over 1-(p+q)}\right)\\[5pt]
q^{min} & = & \min\left({\lambda\over 3}, \mu(1-(p+q)), \mu_i p, \mu_i q\right) = {\lambda\over 3}\min\left(1, {c p\over 1-(p+q)},{c q\over 1-(p+q)} \right)\\[5pt]
\check{q}^{max} & = & 3\lambda_1+3\mu_1 = 3(1+c)\lambda_1 ={\lambda(1+c)\over 1-(p+q)}\\[5pt]
\end{array}
$$
For the resulting bound with $\lambda=1$, $c$ ranging from $2$ to $9$ and for $p,q$ close to 0, the ratio 
of the spectral gap and  (\ref{final-1-2})  in the best case  is of order $10^{-5}$. In this example the bound  (\ref{final-1-2}) is rather rough.
\end{ex}
\begin{ex}
Let $\mathbf{Z}$ be the classical {\em completely symmetrical}  Jackson network with $m$ stations, the  routing matrix $R$  given by $r_{ij}=p<1/(m-1)$ for all $i\neq j$, $r_{0i}=1/m, i,j=1,\ldots,m$, and the 
arrival intensity $\lambda$. Note, that we have  $r_{i0}=1-(m-1)p$ for $i=1,\ldots,m$. The solution of the traffic equation is given by $\lambda_i={1\over m} {\lambda\over 1-(m-1)p}$ for all $i=1,\ldots,m$.
Moreover, assume that $\mu_i=c\lambda_i, c>1$. Then the assumptions of Proposition \ref{gap3} $(i)$ are fulfilled and
$$Gap({\bf Q}({\bf Z})) \le Gap_{ess}:= \left(1-{(m-1)p^2\over 1-(m-2)p}\right) \lambda_i(\sqrt{c}-1)^2 = $$
$$={1\over m} {1+p\over 1-p(m-2)} (\sqrt{c}-1)^2 \lambda .$$
Note that for $p\in(0,1/(m-1))$ we have
$$ {1\over m} (\sqrt{c}-1)^2 \lambda  \leq Gap_{ess}\leq (\sqrt{c}-1)^2 \lambda$$
Let us compare the value of  the upper bound with the  lower bound obtained in  Proposition \ref{gap2}.
Again, the partial balance condition (\ref{detbal1}) holds, and all birth and death processes $\tilde{Q}_i, i=1,\ldots,m$ are equal in distribution. The intensities are
 $\tilde{\lambda}_i=\lambda r_{0i} = \lambda/m$, and $\tilde{\mu}_i=\mu_i r_{i0}=c\lambda /m$.
 We have (similarly as in the previous example)
$$ Gap(\tilde{Q}_i)=\left(\sqrt{\tilde{\mu}_i}-\sqrt{\tilde{\lambda_i}}\right)^2={\lambda\over m}(\sqrt{c}-1)^2,
$$
therefore
$$
Gap({\bf Q}({\bf Z})) \geq {\lambda\over m}(\sqrt{c}-1)^2.
$$
The obtained bound is the  best we can have as a bound which is independent from $p$. The lower  bound  is at most $m$ times smaller than the above given  upper bound of the spectral gap. Moreover, the exact value $Gap({\bf Q}({\bf Z}))$  can be  arbitrarily close to ${\lambda\over m}(\sqrt{c}-1)^2$ for small values of $p$.
\smallskip\par\noindent
Regarding the bound from Proposition \ref{final},  again each $\pi_i, i=1,\ldots,m$ is geometric 
with the  hazard function $h_{\pi_i}(n)=1-\frac{1}{c}$.  We can reduce $1+\bar{d}\hspace{1pt}\bar{b}(2m+1)$  to 1. We need to calculate the following constants
$$
\begin{array}{lll}
|{\bf Q}| & =  & \lambda\left(1+{c\over 1-(m-1)p}\right),\\[5pt]
q^{min} & = & \min\left({\lambda\over m}, \mu_i(1-(m-1)p), \mu_i p\right) = { \lambda\over m}\min\left(1, {c p\over 1-(m-1)p} \right),\\[5pt]
\check{q}^{max} & = & m\lambda_1+m\mu_1 =  {\lambda(1+c)\over 1-(m-1)p}.\\[5pt]
\end{array}
$$
We skip writing the exact formula for the lower bound.  The resulting values with $\lambda=1$, $c$ ranging from $2$ to $9$ and for $p$ close to 0, compared
to the spectral gap, in the best case,  are of order $10^{-5}$, so the bound  (\ref{final-1-2}) is again rather rough.
\end{ex}

\begin{rem}
 Although the bounds obtained from our  Proposition \ref{final} gave rather rough results it is worth stressing that it is possible to compute them for a large class of networks with variable service rates and unreliable nodes. The results possible to obtain via Proposition \ref{gap3} are limited to very special cases of classical networks with constant service intensities. 
 The bounds from Proposition \ref{gap2} are limited to reliable networks and require a kind of partial balance (\ref{detbal1}) (which is fulfilled for example for reversible networks)  but they are applicable to networks with variable service intensities and seem to work quite well.  It is not true in general that the gap for a network is equal to the gap of a bottleneck station in this network. It still remains a lot of research to do in order to provide good computable bounds for networks especially when the service rates are dependent on the queue size and the nodes can be unreliable.
\end{rem}

{\bf Acknowledgement}.
The authors are grateful to the referee for valuable comments which improved this paper.

\end{document}